\documentclass[12pt]{article}
\usepackage{amsfonts}
\usepackage{graphicx}
\usepackage{amsmath}

\input{tcilatex}
\begin{document}

\begin{center}
{\LARGE Growth of Solutions of Complex Differential Equations in a Sector of
the Unit Disc}

\quad

\textbf{Benharrat BELA\"{I}DI}

\textbf{Department of Mathematics, }

\textbf{Laboratory of Pure and Applied Mathematics, }

\textbf{University of Mostaganem (UMAB), B. P. 227 Mostaganem-Algeria}

\textbf{benharrat.belaidi@univ-mosta.dz}

\textbf{\quad }
\end{center}

\noindent \textbf{Abstract.} In this paper, we deal with the growth of
solutions of homogeneous linear complex differential equation by using the
concept of lower [\textit{p,q}]-order and lower [\textit{p,q}]-type in a
sector of the unit disc instead of the whole unit disc, and we obtain
similar results as in the case of the unit disc.\newline
\newline
\textbf{AMS (2010)} : 34M10, 30D35.\newline

\section{Definitions and introduction}

\noindent Throughout this paper, we shall assume that readers are familiar
with the fundamental results and the standard notations of Nevanlinna's
theory in the complex plane and in the unit disc $\Delta =\left\{ z\in 
\mathbb{C}:\left\vert z\right\vert <1\right\} ,$ see $\left[ 5,6,7,9,15,23%
\right] $.

\quad

\noindent \qquad Consider for $k\geq 2$ the complex differential equation 
\begin{equation}
f^{(k)}(z)+A_{k-1}(z)f^{(k-1)}+\cdots +A_{0}(z)f=0,  \tag{1.1}
\end{equation}%
where coefficients $A_{j}$ ($j=0,1,\dots ,k-1$) are analytic functions in
the unit disc $\Delta .$ It is well-known that every solution of $\left(
1.1\right) $ is analytic in $\Delta $, and there are exactly $k$ linearly
independent solutions of equation $\left( 1.1\right) $ (see e.g. $\left[ 7%
\right] $). The theory of complex differential equations in the unit disc
has been developed since 1980's, see $\left[ 13\right] $. In the year 2000,
Heittokangas in $\left[ 7\right] $ firstly investigated the growth and
oscillation theory of equation $\left( 1.1\right) $ when the coefficients $%
A_{j}$ ($j=0,1,\dots ,k-1$) are analytic functions in the unit disc $\Delta $
by introducing the definition of the function spaces. His results also gave
some important tools for further investigations on the theory of meromorphic
solutions of equations $\left( 1.1\right) $. In 1994, Wu $\left[ 17,18\right]
$ used the Nevanlinna theory in an angle to study the order of growth of
solutions of the second-order linear differential equation in an angular
region. Later Xu and Yi $\left[ 22\right] $, Wu $\left[ 19\right] ,$ Wu and
Li $\left[ 20\right] ,$ Zhang $\left[ 24\right] $ generalized some results
of $\left[ 17,18\right] $ to the case of linear higher order differential
equations in angular domains by using the concepts of iterated $p-$order and
the spread relation. Recently, Wu in $\left[ 21\right] $ developed a new
investigation related to linear differential equations with analytic
coefficients in a sector of the unit disc%
\begin{equation*}
\Omega _{\alpha ,\beta }=\left\{ z\in \mathbb{C}:\alpha <\arg z<\beta ,\text{
}\left\vert z\right\vert <1\right\} ,
\end{equation*}%
and obtained some results about the order of growth of solutions of the
differential equation 
\begin{equation}
A_{k}(z)f^{(k)}(z)+A_{k-1}(z)f^{(k-1)}+\cdots +A_{0}(z)f=0,  \tag{1.2}
\end{equation}%
where coefficients $A_{j}$ ($j=0,1,\dots ,k$) are analytic functions in the
sector $\Omega _{\alpha ,\beta }.$ After that, Long in $\left[ 11,12\right]
, $ Zemirni and Belaidi in $\left[ 25\right] $ obtained different results
concerning the growth of solutions of $\left( 1.1\right) $ and $\left(
1.2\right) $ by using the concepts of iterated $p$-order and [\textit{p,q}%
]-order in the sector $\Omega _{\alpha ,\beta }$. In this paper, we continue
to investigate this new problem and study the growth of solutions of
equation $\left( 1.1\right) $ when the coefficients $A_{j}$ ($j=0,1,\dots
,k-1$) are analytic functions of [\textit{p,q}]-order in the sector $\Omega
_{\alpha ,\beta }.$ Before stating our main results, we give some notations
and basic definitions of meromorphic functions in the unit disc $\Delta $
and in a sector $\Omega _{\alpha ,\beta }$ of the unit disc. The order of a
meromorphic function $f$ in $\Delta $ is defined by%
\begin{equation*}
\rho \left( f\right) =\underset{r\rightarrow 1^{-}}{\lim \sup }\frac{\log
T\left( r,f\right) }{\log \frac{1}{1-r}},
\end{equation*}%
where $T(r,f)$ is the Nevanlinna characteristic function of $f$. If $f$ is
analytic function in $\Delta ,$ then 
\begin{equation*}
\rho _{M}\left( f\right) =\underset{r\rightarrow 1^{-}}{\lim \sup }\frac{%
\log \log M\left( r,f\right) }{\log \frac{1}{1-r}},
\end{equation*}%
where $M(r,f)=\max\limits_{\underset{{\large z\in \Delta }}{|z|=r}%
}\left\vert f(z)\right\vert $ is the maximum modulus function.

\quad

\noindent \textbf{Remark 1.1 }The following two statements hold $[15,$ p. $%
205].$

\noindent $\left( \text{a}\right) $ If $f$ is an analytic function in $%
\Delta $, then%
\begin{equation*}
\rho \left( f\right) \leq \rho _{M}\left( f\right) \leq \rho \left( f\right)
+1
\end{equation*}%
\noindent $\left( \text{b}\right) $ There exist analytic functions $f$ in $%
\Delta $ which satisfy $\rho _{M}\left( f\right) \neq \rho \left( f\right) .$
For example, let $\mu >1$ be a constant, and set 
\begin{equation*}
h\left( z\right) =\exp \left\{ \left( 1-z\right) ^{-\mu }\right\} ,
\end{equation*}%
where we choose the principal branch of the logarithm. Then $\rho \left(
h\right) =\mu -1$ and $\rho _{M}\left( h\right) =\mu $, see $[4]$.

\noindent \qquad In contrast, the possibility that occurs in $\left( \text{b}%
\right) $ cannot occur in the whole plane $%
\mathbb{C}
,$ because if $\rho \left( f\right) $ and $\rho _{M}\left( f\right) $ denote
the order of an entire function $f$ in the plane $%
\mathbb{C}
$ (defined by the Nevanlinna characteristic and the maximum modulus,
respectively), then it is well-know that 
\begin{equation*}
\rho \left( f\right) =\underset{r\rightarrow +\infty }{\lim \sup }\frac{\log
T\left( r,f\right) }{\log r}=\rho _{M}\left( f\right) =\underset{%
r\rightarrow +\infty }{\lim \sup }\frac{\log \log M\left( r,f\right) }{\log r%
}.
\end{equation*}%
The meromorphic function $f$ in the unit disc can be divided into the
following three classes:

\noindent (1) bounded type if $T\left( r,f\right) =O\left( 1\right) $ as $%
r\rightarrow 1^{-};$

\noindent (2) rational or non-admissible type if $T\left( r,f\right)
=O\left( \frac{1}{1-r}\right) $ and $f$ does not belong to (1)$;$

\noindent (3) admissible in $\Delta $\ if%
\begin{equation*}
\underset{r\rightarrow 1^{-}}{\lim \sup }\frac{T\left( r,f\right) }{\log 
\frac{1}{1-r}}=\infty .
\end{equation*}%
\textbf{Definition 1.1 }$\left[ 2,3\right] $ Let $p\geq q\geq 1$ be
integers. Let $f$ be a meromorphic function in $\Delta ,$ the [\textit{p,q}%
]-order \ of $f$ is defined by 
\begin{equation*}
\rho _{\lbrack p,q]}\left( f\right) =\underset{r\rightarrow 1^{-}}{\lim \sup 
}\frac{\log _{p}^{+}T\left( r,f\right) }{\log _{q}\frac{1}{1-r}},
\end{equation*}%
where $\log _{1}^{+}r:=\log ^{+}r=\max \left( 0,\log r\right) ,$ $\log
_{p+1}^{+}r:=\log ^{+}\left( \log _{p}^{+}r\right) ,$ $p\in \mathbb{N}$. For
an analytic function $f$ in $\Delta ,$ we also define 
\begin{equation*}
\rho _{M,[p,q]}\left( f\right) =\underset{r\rightarrow 1^{-}}{\lim \sup }%
\frac{\log _{p+1}^{+}M\left( r,f\right) }{\log _{q}\frac{1}{1-r}}.
\end{equation*}%
It is easy to see that $0\leq \rho _{\lbrack p,q]}\left( f\right) \leq
+\infty .$ If $f$ is non-admissible$,$ then $\rho _{\lbrack p,q]}\left(
f\right) =0$ for any $p\geq q\geq 1.$ By Definition 1.1, $\rho _{\lbrack
1,1]}\left( f\right) =\rho \left( f\right) $ is the order of $f$ in $\Delta
, $ $\rho _{\lbrack 2,1]}\left( f\right) =\rho _{2}\left( f\right) $ is the
hyper-order of $f$ in $\Delta $ and $\rho _{\lbrack p,1]}\left( f\right)
=\rho _{p}\left( f\right) $ is the $p$-iterated order of $f$ in $\Delta .$

\quad

\noindent \textbf{Proposition 1.1} $\left[ 2\right] $ Let $p\geq q\geq 1$ be
integers, and let $f$ be an analytic function in $\Delta $ of [\textit{p,q}%
]-order. The following two statements hold :\newline
$(i)$ If $p=q$, then%
\begin{equation*}
\rho _{\lbrack p,q]}\left( f\right) \leq \rho _{M,[p,q]}\left( f\right) \leq
\rho _{\lbrack p,q]}\left( f\right) +1.
\end{equation*}%
$(ii)$ If $p>q,$ then%
\begin{equation*}
\rho _{\lbrack p,q]}\left( f\right) =\rho _{M,[p,q]}\left( f\right) .
\end{equation*}%
\textbf{Proposition 1.2 }$\left[ 8\right] $ Let $p\geq q\geq 1$ be integers,
and let $f$ be an analytic function in $\Delta $ of [\textit{p,q}]-order.
The following two statements hold :\newline
$(i)$ If $p=q$, then%
\begin{equation*}
\mu _{\lbrack p,q]}\left( f\right) \leq \mu _{M,[p,q]}\left( f\right) \leq
\mu _{\lbrack p,q]}\left( f\right) +1.
\end{equation*}%
$(ii)$ If $p>q,$ then%
\begin{equation*}
\mu _{\lbrack p,q]}\left( f\right) =\mu _{M,[p,q]}\left( f\right) .
\end{equation*}

\noindent \qquad In what follows, we give some notations and definitions of
a meromorphic function in a sector in unit disc. Throughout this paper, $%
\Omega $ usually denotes the sector $\Omega _{\alpha ,\beta }$ ($0\leq
\alpha <\beta \leq 2\pi $) of the unit disc, and for any given $\varepsilon
\in \left( 0,\frac{\beta -\alpha }{2}\right) ,$ $\Omega _{\varepsilon }$
denotes the sector%
\begin{equation*}
\Omega _{\alpha ,\beta ,\varepsilon }=\left\{ z\in \mathbb{C}:\alpha
+\varepsilon <\arg z<\beta -\varepsilon ,\text{ }\left\vert z\right\vert
<1\right\} .
\end{equation*}%
In $\left[ 21\right] $, Wu has used the Ahlfors-Shimizu characteristic
function to measure the order of growth of a meromorphic function $f$ in $%
\Omega .$ We recall the definition of the Ahlfors-Shimizu characteristic
function, see $\left[ 5,6\right] $. Let $f$ be a meromorphic function in $%
\Omega ,$ set%
\begin{equation*}
\Omega (r)=\Omega \cap \left\{ z\in \mathbb{C}:0<\left\vert z\right\vert
<r<1\right\}
\end{equation*}%
\begin{equation*}
=\left\{ z\in \mathbb{C}:\alpha <\arg z<\beta ,\text{ }0<\left\vert
z\right\vert <r<1\right\} .
\end{equation*}%
Then, the Ahlfors-Shimizu characteristic function is defined by%
\begin{equation*}
T_{0}\left( r,\Omega ,f\right) =\int_{0}^{r}\frac{S\left( t,\Omega ,f\right) 
}{t}dt,
\end{equation*}%
where%
\begin{equation*}
S\left( r,\Omega ,f\right) =\frac{1}{\pi }\iint\limits_{\Omega (r)}\left( 
\frac{\left\vert f^{\prime }\left( z\right) \right\vert }{1+\left\vert
f\left( z\right) \right\vert ^{2}}\right) ^{2}d\sigma ,\text{ }z=re^{i\theta
},\text{ }d\sigma =rdrd\theta .
\end{equation*}%
It follows by Hayman $\left[ 6\right] $, Goldberg and Ostrovskii $\left[ 5%
\right] $ that 
\begin{equation*}
T_{0}\left( r,\mathbb{C},f\right) =T\left( r,f\right) +O\left( 1\right) ,%
\text{ }0<r<1.
\end{equation*}%
The meromorphic function $f$ in a sector $\Omega $ of the unit disc can be
divided into the following three classes:

\noindent (1) bounded type if $T_{0}\left( r,\Omega ,f\right) =O\left(
1\right) $ as $r\rightarrow 1^{-};$

\noindent (2) rational or non-admissible type if $T_{0}\left( r,\Omega
,f\right) =O\left( \frac{1}{1-r}\right) $ and $f$ does not belong to (1);

\noindent (3) admissible in $\Omega $ if%
\begin{equation*}
\underset{r\rightarrow 1^{-}}{\lim \sup }\frac{T_{0}\left( r,\Omega
,f\right) }{\log \frac{1}{1-r}}=\infty .
\end{equation*}

\noindent \qquad Now, we introduce the concept of [\textit{p,q}]-order and [%
\textit{p,q}]-type of meromorphic functions in a sector $\Omega $.

\quad

\noindent \textbf{Definition 1.2 }$\left[ 12,25\right] $ Let $p\geq q\geq 1$
be integers. Let $f$ be a meromorphic function in $\Omega ,$ the [\textit{p,q%
}]-order\ of $f$ is defined by 
\begin{equation*}
\rho _{\lbrack p,q],\Omega }\left( f\right) =\underset{r\rightarrow 1^{-}}{%
\lim \sup }\frac{\log _{p}^{+}T_{0}\left( r,\Omega ,f\right) }{\log _{q}%
\frac{1}{1-r}}.
\end{equation*}%
It is clear that $0\leq \rho _{\lbrack p,q],\Omega }\left( f\right) \leq
+\infty .$ If $f$ is non-admissible in $\Omega ,$ then $\rho _{\lbrack
p,q],\Omega }\left( f\right) =0.$ By Definition 1.2, $\rho _{\lbrack
1,1],\Omega }\left( f\right) =\rho _{\Omega }\left( f\right) $ is the order
of $f$ in $\Omega ,$ see $\left[ 21\right] $, $\rho _{\lbrack p,1],\Omega
}\left( f\right) =\rho _{p,\Omega }\left( f\right) $ is the iterated $p$%
-order of $f$ in $\Omega ,$ see $\left[ 11,24\right] $.

\quad

\noindent \textbf{Definition 1.3 }$\left[ 25\right] $ Let $p\geq q\geq 1$ be
integers and $f$ be a meromorphic function in $\Omega $ with [\textit{p,q}%
]-order $0<\rho _{\lbrack p,q],\Omega }\left( f\right) <+\infty .$ Then, the
[\textit{p,q}]-type\ of $f$ is defined by 
\begin{equation*}
\tau _{\lbrack p,q],\Omega }\left( f\right) =\underset{r\rightarrow 1^{-}}{%
\lim \sup }\frac{\log _{p-1}^{+}T_{0}\left( r,\Omega ,f\right) }{\left( \log
_{q-1}\frac{1}{1-r}\right) ^{\rho _{\lbrack p,q],\Omega }\left( f\right) }}.
\end{equation*}

\noindent \qquad Now, we introduce the concept of lower [\textit{p,q}]-order
and lower [\textit{p,q}]-type of a meromorphic function in a sector $\Omega $%
.

\quad

\noindent \textbf{Definition 1.4} Let $p\geq q\geq 1$ be integers. Let $f$
be a meromorphic function in $\Omega ,$ the lower [\textit{p,q}]-order\ of $%
f $ is defined by 
\begin{equation*}
\mu _{\lbrack p,q],\Omega }\left( f\right) =\underset{r\rightarrow 1^{-}}{%
\lim \inf }\frac{\log _{p}^{+}T_{0}\left( r,\Omega ,f\right) }{\log _{q}%
\frac{1}{1-r}}.
\end{equation*}%
It is clear that $0\leq \mu _{\lbrack p,q],\Omega }\left( f\right) \leq
+\infty .$ If $f$ is non-admissible in $\Omega ,$ then $\mu _{\lbrack
p,q],\Omega }\left( f\right) =0.$ By Definition 1.4, $\mu _{\lbrack
1,1],\Omega }\left( f\right) =\mu _{\Omega }\left( f\right) $ is the lower
order of $f$ in $\Omega $ and $\mu _{\lbrack p,1],\Omega }\left( f\right)
=\mu _{p,\Omega }\left( f\right) $ is the lower iterated $p$-order of $f$ in 
$\Omega .$

\quad

\noindent \textbf{Definition 1.5} Let $p\geq q\geq 1$ be integers and $f$ be
a meromorphic function in $\Omega $ with lower [\textit{p,q}]-order $0<\mu
_{\lbrack p,q],\Omega }\left( f\right) <+\infty .$ Then, the lower [\textit{%
p,q}]-type\ of $f$ is defined by 
\begin{equation*}
\underline{\tau }_{[p,q],\Omega }\left( f\right) =\underset{r\rightarrow
1^{-}}{\lim \inf }\frac{\log _{p-1}^{+}T_{0}\left( r,\Omega ,f\right) }{%
\left( \log _{q-1}\frac{1}{1-r}\right) ^{\mu _{\lbrack p,q],\Omega }\left(
f\right) }}.
\end{equation*}

\section{\textbf{Main results}}

\noindent \qquad Several authors $\left[ 2,3,8,10,16\right] $ have
investigated the growth of solutions of the equation $\left( 1.1\right) $ by
using the concepts of [\textit{p,q}]-order in the unit disc $\Delta $. In $%
\left[ 12\right] $, Long has studied the growth of solutions of the equation 
$\left( 1.2\right) $ in a sector of the unit disc with analytic coefficients
of finite [\textit{p,q}]-order, and has obtained the following results.

\quad

\noindent \textbf{Theorem A }$\left[ 12\right] $ \textit{Let }$p\geq q\geq 1$%
\textit{\ be integers and }$\varepsilon \in \left( 0,\frac{\beta -\alpha }{2}%
\right) $\textit{. Let }$E$\textit{\ be a set of complex numbers satisfying }%
$\overline{\mathrm{dens}}\left\{ \left\vert z\right\vert =r:z\in E\subset
\Omega \right\} >0,$\textit{\ and let }$A_{0}(z),A_{1}(z),\dots ,A_{k}(z)$%
\textit{\ be analytic functions in }$\Omega $\textit{\ such that for some
real constants satisfying }$0\leq \gamma <\lambda $\textit{, we have}%
\begin{equation*}
T_{0}\left( r,\Omega _{\varepsilon },A_{0}(z)\right) \geq \exp _{p}\left\{
\lambda \log _{q}\left( \frac{1}{1-\left\vert z\right\vert }\right) \right\}
,
\end{equation*}%
\begin{equation*}
T_{0}\left( r,\Omega ,A_{j}(z)\right) \leq \exp _{p}\left\{ \gamma \log
_{q}\left( \frac{1}{1-\left\vert z\right\vert }\right) \right\} ,\text{ }%
j=1,2,\dots ,k
\end{equation*}%
\textit{as }$\left\vert z\right\vert =r\rightarrow 1^{-}$\textit{\ for }$%
z\in E.$\textit{\ Then every nontrivial solution }$f$\textit{\ of }$\left(
1.2\right) $\textit{\ satisfies }$\rho _{\lbrack p,q],\Omega }\left(
f\right) =+\infty $\ \textit{and }%
\begin{equation*}
\rho _{\lbrack p+1,q],\Omega }\left( f\right) \geq \lambda .
\end{equation*}%
\textbf{Theorem B }$\left[ 12\right] $\textbf{\ }\textit{Let }$p\geq q\geq 1$%
\textit{\ be integers and }$\varepsilon \in \left( 0,\frac{\beta -\alpha }{2}%
\right) $\textit{. Let }$A_{0}(z),A_{1}(z),\dots ,A_{k}(z)$\textit{\ be
analytic functions in }$\Omega .$\textit{\ If}%
\begin{equation*}
\max_{1\leq j\leq k}\left\{ \rho _{\lbrack p,q],\Omega }\left( A_{j}\right)
\right\} <\rho _{\lbrack p,q],\Omega _{\varepsilon }}\left( A_{0}\right) ,
\end{equation*}%
\textit{then every nontrivial solution of }$\left( 1.2\right) $\textit{\
satisfies }%
\begin{equation*}
\rho _{\lbrack p+1,q],\Omega }\left( f\right) \geq \rho _{\lbrack
p,q],\Omega _{\varepsilon }}\left( A_{0}\right) .
\end{equation*}%
\textbf{Remark 2.1.} In Theorems A-B, we note that if $A_{k}(z)=1$, then all
the solutions of equation $\left( 1.2\right) $ are analytic functions. But
if $A_{k}(z)$ is a non constant analytic function, then obviously the
solution $f$ of the equation $\left( 1.2\right) $ can be meromorphic
function.\textbf{\ }The\textbf{\ }hypotheses of the Theorems A-B do not
provide that a solution is meromorphic in $\Omega $, so it is a priori
assumed that $f$ is meromorphic.

\quad

\noindent Very recently, Zemirni and Bela\"{\i}di $\left[ 25\right] $ have
continued the study of the growth of solutions of the equation $\left(
1.1\right) $ instead of the equation $\left( 1.2\right) $ in a sector of the
unit disc with analytic coefficients of finite [\textit{p,q}]-order, and
have got the following results.

\quad

\noindent \textbf{Theorem C }$\left[ 25\right] $ \textit{Let }$p\geq q\geq 1$%
\textit{\ be integers and }$\varepsilon \in \left( 0,\frac{\beta -\alpha }{2}%
\right) $\textit{. Let }$A_{0}(z),A_{1}(z),\newline
\dots ,$ $A_{k-1}(z)$\textit{\ be analytic functions in }$\Omega .$\textit{\
If}%
\begin{equation*}
\max_{1\leq j\leq k-1}\left\{ \rho _{\lbrack p,q],\Omega }\left(
A_{j}\right) \right\} <\rho _{\lbrack p,q],\Omega _{\varepsilon }}\left(
A_{0}\right) ,
\end{equation*}%
\textit{then every nontrivial solution of }$\left( 1.1\right) $\textit{\
satisfies }$\rho _{\lbrack p,q],\Omega }\left( f\right) =+\infty $\textit{\
and }%
\begin{equation*}
\rho _{\lbrack p,q],\Omega _{\varepsilon }}\left( A_{0}\right) \leq \rho
_{\lbrack p+1,q],\Omega }\left( f\right) \text{, }\rho _{\lbrack
p+1,q],\Omega _{\varepsilon }}\left( f\right) \leq \rho _{\lbrack
p,q],\Omega }\left( A_{0}\right) +1.
\end{equation*}%
\textit{Furthermore, if }$p>q,$\textit{\ then}%
\begin{equation*}
\rho _{\lbrack p,q],\Omega _{\varepsilon }}\left( A_{0}\right) \leq \rho
_{\lbrack p+1,q],\Omega }\left( f\right) \text{,\ }\rho _{\lbrack
p+1,q],\Omega _{\varepsilon }}\left( f\right) \leq \rho _{\lbrack
p,q],\Omega }\left( A_{0}\right) .
\end{equation*}%
\textbf{Theorem D} $\left[ 25\right] $ \textit{Let }$p\geq q\geq 1$\textit{\
be integers and }$\varepsilon \in \left( 0,\frac{\beta -\alpha }{2}\right) $%
\textit{. Let }$A_{0}(z),A_{1}(z),\newline
\dots ,$ $A_{k-1}(z)$\textit{\ be analytic functions in }$\Omega .$\textit{\
Suppose that }%
\begin{equation*}
\max_{1\leq j\leq k-1}\left\{ \rho _{\lbrack p,q],\Omega }\left(
A_{j}\right) \right\} \leq \rho _{\lbrack p,q],\Omega _{\varepsilon }}\left(
A_{0}\right) =\rho \text{ }\left( 0<\rho <+\infty \right)
\end{equation*}%
\textit{and }%
\begin{equation*}
\max_{1\leq j\leq k-1}\{\tau _{\lbrack p,q],\Omega }\left( A_{j}\right)
:\rho _{\lbrack p,q],\Omega }\left( A_{j}\right) =\rho _{\lbrack p,q],\Omega
_{\varepsilon }}\left( A_{0}\right) \}
\end{equation*}%
\begin{equation*}
<\tau _{\lbrack p,q],\Omega _{\varepsilon }}\left( A_{0}\right) =\tau \text{ 
}(0<\tau <+\infty ).
\end{equation*}%
\textit{Then every nontrivial solution of }$\left( 1.1\right) $\textit{\
satisfies }$\rho _{\lbrack p,q],\Omega }\left( f\right) =+\infty $\textit{\
and }%
\begin{equation*}
\rho _{\lbrack p,q],\Omega _{\varepsilon }}\left( A_{0}\right) \leq \rho
_{\lbrack p+1,q],\Omega }\left( f\right) ,\text{ }\rho _{\lbrack
p+1,q],\Omega _{\varepsilon }}\left( f\right) \leq \rho _{\lbrack
p,q],\Omega }\left( A_{0}\right) +1.
\end{equation*}%
\textit{Furthermore, if }$p>q,$\textit{\ then}%
\begin{equation*}
\rho _{\lbrack p,q],\Omega _{\varepsilon }}\left( A_{0}\right) \leq \rho
_{\lbrack p+1,q],\Omega }\left( f\right) ,\text{ \ }\rho _{\lbrack
p+1,q],\Omega _{\varepsilon }}\left( f\right) \leq \rho _{\lbrack
p,q],\Omega }\left( A_{0}\right) .
\end{equation*}

\noindent \qquad Thus, the following questions arise naturally: (i) Whether
the results similar to Theorem C can be obtained in\textit{\ }$\Omega $ if $%
A_{0}\left( z\right) $ to dominate other coefficients in the sense of lower [%
\textit{p,q}]-order?

\noindent (ii) If we use the lower [\textit{p,q}]-type of $A_{0}\left(
z\right) $ to dominate other coefficients, what can be said about $\mu
_{\lbrack p+1,q],\Omega }\left( f\right) $ similar to Theorem D? In this
paper, we give some answers to the above questions. In fact, by using the
concept of lower [\textit{p,q}]-type, we obtain some results which indicate
growth estimate of every non-trivial analytic solution of equation $\left(
1.1\right) $ by the growth estimate of the coefficient $A_{0}\left( z\right)
.$We mainly obtain the following results.

\quad

\noindent \textbf{Theorem 2.1} \textit{Let }$p\geq q\geq 1$\textit{\ be
integers and }$\varepsilon \in \left( 0,\frac{\beta -\alpha }{2}\right) $%
\textit{. Let }$A_{0}(z),A_{1}(z),\newline
\dots ,$ $A_{k-1}(z)$\textit{\ be analytic functions in }$\Omega .$\textit{\
If}%
\begin{equation*}
\max_{1\leq j\leq k-1}\left\{ \rho _{\lbrack p,q],\Omega }\left(
A_{j}\right) \right\} <\mu _{\lbrack p,q],\Omega _{\varepsilon }}\left(
A_{0}\right) ,
\end{equation*}%
\textit{then every nontrivial solution of }$\left( 1.1\right) $\textit{\
satisfies }$\rho _{\lbrack p,q],\Omega }\left( f\right) =\mu _{\lbrack
p,q],\Omega }\left( f\right) =+\infty ,$\textit{\ }%
\begin{equation*}
\mu _{\lbrack p,q],\Omega _{\varepsilon }}\left( A_{0}\right) \leq \mu
_{\lbrack p+1,q],\Omega }\left( f\right) \leq \rho _{\lbrack p+1,q],\Omega
}\left( f\right) \text{ }
\end{equation*}%
\textit{and }%
\begin{equation*}
\mu _{\lbrack p+1,q],\Omega _{\varepsilon }}\left( f\right) \leq \mu
_{\lbrack p,q],\Omega }\left( A_{0}\right) +1.
\end{equation*}%
\textit{Furthermore, if }$p>q,$\textit{\ then}%
\begin{equation*}
\mu _{\lbrack p,q],\Omega _{\varepsilon }}\left( A_{0}\right) \leq \mu
_{\lbrack p+1,q],\Omega }\left( f\right) \leq \rho _{\lbrack p+1,q],\Omega
}\left( f\right)
\end{equation*}%
\textit{and }%
\begin{equation*}
\text{ }\mu _{\lbrack p+1,q],\Omega _{\varepsilon }}\left( f\right) \leq \mu
_{\lbrack p,q],\Omega }\left( A_{0}\right) .
\end{equation*}%
\textbf{Remark 2.2} The Theorem 2.1 is similar to Theorem 2.2 $\left( \text{i%
}\right) $ in $\left[ 16\right] $ in the unit disc $\Delta $.$\newline
$

~

\noindent \textbf{Corollary 2.1} \textit{Let }$p\geq q\geq 1$\textit{\ be
integers and }$\varepsilon \in \left( 0,\frac{\beta -\alpha }{2}\right) $%
\textit{. Let }$A_{0}(z),A_{1}(z),\newline
\dots ,$ $A_{k-1}(z)$\textit{\ be analytic functions in }$\Omega .$\textit{\
If}%
\begin{equation*}
\max_{1\leq j\leq k-1}\left\{ \rho _{\lbrack p,q],\Omega }\left(
A_{j}\right) \right\} <\mu _{\lbrack p,q],\Omega _{\varepsilon }}\left(
A_{0}\right) =\rho _{\lbrack p,q],\Omega _{\varepsilon }}\left( A_{0}\right)
,
\end{equation*}%
\textit{then every nontrivial solution of }$\left( 1.1\right) $\textit{\
satisfies }$\rho _{\lbrack p,q],\Omega }\left( f\right) =\mu _{\lbrack
p,q],\Omega }\left( f\right) =+\infty $\textit{\ and }%
\begin{equation*}
\mu _{\lbrack p,q],\Omega _{\varepsilon }}\left( A_{0}\right) \leq \mu
_{\lbrack p+1,q],\Omega }\left( f\right) \leq \rho _{\lbrack p+1,q],\Omega
}\left( f\right) ,\text{ }
\end{equation*}%
\textit{\ }%
\begin{equation*}
\mu _{\lbrack p+1,q],\Omega _{\varepsilon }}\left( f\right) \leq \rho
_{\lbrack p+1,q],\Omega _{\varepsilon }}\left( f\right) \leq \mu _{\lbrack
p,q],\Omega }\left( A_{0}\right) +1.
\end{equation*}%
\textit{Furthermore, if }$p>q,$\textit{\ then}%
\begin{equation*}
\mu _{\lbrack p,q],\Omega _{\varepsilon }}\left( A_{0}\right) \leq \mu
_{\lbrack p+1,q],\Omega }\left( f\right) \leq \rho _{\lbrack p+1,q],\Omega
}\left( f\right)
\end{equation*}%
\textit{and }%
\begin{equation*}
\text{ }\mu _{\lbrack p+1,q],\Omega _{\varepsilon }}\left( f\right) \leq
\rho _{\lbrack p+1,q],\Omega _{\varepsilon }}\left( f\right) \leq \mu
_{\lbrack p,q],\Omega }\left( A_{0}\right) .
\end{equation*}%
\textbf{Theorem 2.2} \textit{Let }$p\geq q\geq 1$\textit{\ be integers and }$%
\varepsilon \in \left( 0,\frac{\beta -\alpha }{2}\right) $\textit{. Let }$%
A_{0}(z),A_{1}(z),\newline
\dots ,$ $A_{k-1}(z)$\textit{\ be analytic functions in }$\Omega $\textit{\
such that }$0<\mu =\mu _{\lbrack p,q],\Omega _{\varepsilon }}\left(
A_{0}\right) \leq \rho _{\lbrack p,q],\Omega _{\varepsilon }}\left(
A_{0}\right) <+\infty .$\textit{\ Suppose that }%
\begin{equation*}
\max_{1\leq j\leq k-1}\left\{ \rho _{\lbrack p,q],\Omega }\left(
A_{j}\right) \right\} \leq \mu _{\lbrack p,q],\Omega _{\varepsilon }}\left(
A_{0}\right)
\end{equation*}%
\textit{and }%
\begin{equation*}
\max_{1\leq j\leq k-1}\{\tau _{\lbrack p,q],\Omega }\left( A_{j}\right)
:\rho _{\lbrack p,q],\Omega }\left( A_{j}\right) =\mu _{\lbrack p,q],\Omega
_{\varepsilon }}\left( A_{0}\right) \}<\underline{\tau }_{[p,q],\Omega
_{\varepsilon }}\left( A_{0}\right) <+\infty .
\end{equation*}%
\textit{Then every nontrivial solution of }$\left( 1.1\right) $\textit{\
satisfies }$\rho _{\lbrack p,q],\Omega }\left( f\right) =\mu _{\lbrack
p,q],\Omega }\left( f\right) =+\infty $\textit{\ and }%
\begin{equation*}
\mu _{\lbrack p,q],\Omega _{\varepsilon }}\left( A_{0}\right) \leq \mu
_{\lbrack p+1,q],\Omega }\left( f\right) \leq \rho _{\lbrack p+1,q],\Omega
}\left( f\right) ,
\end{equation*}%
\begin{equation*}
\text{ }\mu _{\lbrack p+1,q],\Omega _{\varepsilon }}\left( f\right) \leq \mu
_{\lbrack p,q],\Omega }\left( A_{0}\right) +1.
\end{equation*}%
\textit{Furthermore, if }$p>q,$\textit{\ then}%
\begin{equation*}
\mu _{\lbrack p,q],\Omega _{\varepsilon }}\left( A_{0}\right) \leq \mu
_{\lbrack p+1,q],\Omega }\left( f\right) \leq \rho _{\lbrack p+1,q],\Omega
}\left( f\right) \text{ }
\end{equation*}%
\textit{and}%
\begin{equation*}
\mu _{\lbrack p+1,q],\Omega _{\varepsilon }}\left( f\right) \leq \mu
_{\lbrack p,q],\Omega }\left( A_{0}\right) .
\end{equation*}%
\textbf{Remark 2.3} The Theorem 2.2 is similar to Theorem 2.1 in $\left[ 8%
\right] $ in the unit disc $\Delta $.

\quad

\noindent \textbf{Remark 2.4} We note that in Theorems 2.1 and 2.2, the
growth estimate of the solution $f$ is expressed by the growth estimate of
dominant coefficient $A_{0}$ in the terms of lower [\textit{p,q}]-order on
both sides.

\section{\textbf{Auxiliary lemmas}}

\noindent \textbf{Lemma 3.1} $\left[ 14\right] $ \textit{Let}%
\begin{equation}
u(z)=\frac{\left( ze^{-i\theta _{0}}\right) ^{\pi /\delta }+2\left(
ze^{-i\theta _{0}}\right) ^{\pi /\left( 2\delta \right) }-1}{\left(
ze^{-i\theta _{0}}\right) ^{\pi /\delta }-2\left( ze^{-i\theta _{0}}\right)
^{\pi /\left( 2\delta \right) }-1},  \tag{3.1}
\end{equation}%
\textit{where }$0\leq \theta _{0}=\frac{\alpha +\beta }{2}<2\pi ,$\textit{\ }%
$0<\delta =\frac{\beta -\alpha }{2}<\pi .$\textit{\ Then }$u(z)$\textit{\ is
a conformal map of angular domain }$\Omega ,$\textit{\ }$(0<\beta -\alpha
<2\pi )$\textit{\ onto the unit disc }$\Delta .$\textit{\ Moreover, for any
positive number }$\varepsilon $\textit{\ satisfying }$0<\varepsilon <\delta
, $\textit{\ the transformation }$\left( 3.1\right) $\textit{\ satisfies}%
\begin{equation*}
u\left( \left\{ z:\frac{1}{2}<\left\vert z\right\vert <r\right\} \cap
\left\{ z:\left\vert \arg z-\theta _{0}\right\vert <\delta -\varepsilon
\right\} \right)
\end{equation*}%
\begin{equation*}
\subset \left\{ u:\left\vert u\right\vert <1-\frac{\varepsilon }{2^{\frac{%
\pi }{2\delta }+1}\delta }\left( 1-r\right) \right\} ,
\end{equation*}%
\begin{equation*}
u^{-1}\left( \left\{ u:\left\vert u\right\vert <\varrho \right\} \right)
\subset \left( \left\{ z:\left\vert z\right\vert <1-\frac{\delta }{8\pi }%
(1-\varrho )\right\} \cap \left\{ z:\left\vert \arg z-\theta _{0}\right\vert
<\delta \right\} \right) ,
\end{equation*}%
\textit{where }$\varrho <1$\textit{\ is a constant. The inverse
transformation of }$\left( 3.1\right) $\textit{\ is}%
\begin{equation}
z(u)=e^{i\theta _{0}}\left( \frac{-(1+u)+\sqrt{2(1+u^{2})}}{1-u}\right) ^{%
\frac{2\delta }{\pi }}.  \tag{3.2}
\end{equation}%
\textbf{Lemma 3.2 }$\left[ 21\right] $\textbf{\ }\textit{Let }$f$\textit{\
be a meromorphic function in }$\Omega ,$\textit{\ where }$0<\beta -\alpha
<2\pi .$\textit{\ For any given }$\varepsilon \in \left( 0,\frac{\beta
-\alpha }{2}\right) ,$\textit{\ set }$\delta =\frac{\beta -\alpha }{2}$%
\textit{\ and }$b=\frac{\varepsilon }{2^{\pi /(2\delta )+1}\delta }.$\textit{%
\ Then the following statements hold}%
\begin{equation}
T_{0}\left( \varrho ,\mathbb{C},f\left( z\left( u\right) \right) \right)
\leq \frac{16\pi }{\delta }T_{0}\left( 1-\frac{\delta }{8\pi }\left(
1-\varrho \right) ,\Omega ,f(z)\right) +O(1),  \tag{3.3}
\end{equation}%
\begin{equation}
T_{0}\left( r,\Omega _{\varepsilon },f(z)\right) \leq \frac{2}{b}T_{0}\left(
1-b\left( 1-r\right) ,\mathbb{C},f\left( z(u)\right) \right) +O(1), 
\tag{3.4}
\end{equation}%
\textit{where }$z(u)$\textit{\ is the inverse transformation of }$\left(
3.1\right) .$

\quad

\noindent \textbf{Remark 3.1} \textit{By applying the formula }$T\left(
r,f\right) =T_{0}\left( r,\mathbb{C},f\right) +O(1)$ $\left( 0<r<1\right) ,$%
\textit{\ Lemma 3.2, the definition of }[\textit{p,q}]\textit{-order and
lower }[\textit{p,q}]\textit{-order, we immediately obtain that}%
\begin{equation*}
\rho _{\lbrack p,q],\Omega _{\varepsilon }}\left( f\left( z\right) \right)
\leq \rho _{\lbrack p,q]}\left( f\left( z\left( u\right) \right) \right)
\leq \rho _{\lbrack p,q],\Omega }\left( f\left( z\right) \right)
\end{equation*}%
\textit{and}%
\begin{equation*}
\mu _{\lbrack p,q],\Omega _{\varepsilon }}\left( f\left( z\right) \right)
\leq \mu _{\lbrack p,q]}\left( f\left( z\left( u\right) \right) \right) \leq
\mu _{\lbrack p,q],\Omega }\left( f\left( z\right) \right) .
\end{equation*}%
\textbf{Lemma 3.3 }$\left[ 21\right] $ \textbf{\ }\textit{Let }$f$\textit{\
be a meromorphic function in }$\Omega ,$\textit{\ where }$0<\beta -\alpha
<2\pi $\textit{\ and }$z(u)$\textit{\ be the inverse transformation of }$%
\left( 3.1\right) .$\textit{\ Set }$F(u)=f\left( z\left( u\right) \right) ,$%
\textit{\ }$\psi \left( u\right) =f^{(\ell )}\left( z\left( u\right) \right)
,$\textit{\ then}%
\begin{equation}
\psi \left( u\right) =\sum_{j=1}^{\ell }\alpha _{j}F^{(j)}(u),  \tag{3.5}
\end{equation}%
\textit{where the coefficients }$\alpha _{j}$\textit{\ are polynomials (with
numerical coefficients) in the variables }$V(u)\left( =\frac{1}{z^{\prime
}(u)}\right) ,V^{\prime }(u),V^{\prime \prime }(u),\dots .$\textit{\
Moreover, we have}%
\begin{equation}
T\left( \varrho ,\alpha _{j}\right) =O\left( \log \frac{1}{1-\varrho }%
\right) ,\text{ \ }j=1,2,\dots ,\ell .  \tag{3.6}
\end{equation}%
For the convenience of the readers, we give the statement and the proof of
Lemma 3.4 $\left[ 25,\text{ Lemma 3.4}\right] $ with more precisions.

\quad

\noindent \textbf{Lemma 3.4 }\textit{Suppose }$f\not\equiv 0$\textit{\ is a
solution of }$\left( 1.1\right) $\textit{\ in }$\Omega .$\textit{\ Then }$%
F(u)=f\left( z\left( u\right) \right) $\textit{\ is a solution of }%
\begin{equation}
F^{(k)}(u)+B_{k-1}(u)F^{(k-1)}(u)+\cdots +B_{0}(u)F(u)=0  \tag{3.7}
\end{equation}%
\textit{in }$\Delta ,$\textit{\ where }%
\begin{equation}
B_{0}(u)=\dfrac{1}{\alpha _{k}}A_{0}\left( z\left( u\right) \right) 
\tag{3.8}
\end{equation}%
\textit{\ and for }$j=1,2,\dots ,k-1$%
\begin{equation}
B_{j}(u)=\frac{\alpha _{j}}{\alpha _{k}}+\frac{\alpha _{j}}{\alpha _{k}}%
\sum_{n=j}^{k-1}A_{n}\left( z\left( u\right) \right) .  \tag{3.9}
\end{equation}%
\textit{\ Consequently},%
\begin{equation}
T\left( \varrho ,B_{0}\right) \leq T\left( r,A_{0}\left( z\left( u\right)
\right) \right) +O\left( \log \frac{1}{1-\varrho }\right)  \tag{3.10}
\end{equation}%
\textit{and}%
\begin{equation}
T\left( \varrho ,B_{j}\right) \leq \sum_{n=j}^{k-1}T\left( r,A_{n}\left(
z\left( u\right) \right) \right) +O\left( \log \frac{1}{1-\varrho }\right) .
\tag{3.11}
\end{equation}%
\textit{Proof}. Suppose that $f\not\equiv 0$ is a solution of $\left(
1.1\right) $ in the sector $\Omega .$ By using Lemma 3.3, we have%
\begin{equation*}
f^{(k)}(z\left( u\right) )+\sum_{n=1}^{k-1}A_{n}\left( z\left( u\right)
\right) f^{(n)}(z\left( u\right) )+A_{0}\left( z\left( u\right) \right)
f\left( z\left( u\right) \right)
\end{equation*}%
\begin{equation*}
=\sum_{j=1}^{k}\alpha _{j}F^{(j)}(u)+\sum_{n=1}^{k-1}A_{n}\left( z\left(
u\right) \right) \sum_{j=1}^{n}\alpha _{j}F^{(j)}(u)+A_{0}\left( z\left(
u\right) \right) f\left( z\left( u\right) \right)
\end{equation*}%
\begin{equation*}
=\sum_{j=1}^{k}\alpha _{j}F^{(j)}(u)+\sum_{j=1}^{k-1}\left( \alpha
_{j}\sum_{n=j}^{k-1}A_{n}\left( z\left( u\right) \right) \right)
F^{(j)}(u)+A_{0}\left( z\left( u\right) \right) f\left( z\left( u\right)
\right)
\end{equation*}%
\begin{equation*}
=\alpha _{k}F^{(k)}(u)+\sum_{j=1}^{k-1}\left( \alpha
_{j}\sum_{n=j}^{k-1}A_{n}\left( z\left( u\right) \right) +\alpha _{j}\right)
F^{(j)}(u)+A_{0}\left( z\left( u\right) \right) F\left( u\right) .
\end{equation*}%
It follows that $F(u)=f\left( z\left( u\right) \right) $\textit{\ }is a
solution of\textit{\ }%
\begin{equation*}
F^{(k)}(u)+B_{k-1}(u)F^{(k-1)}(u)+\cdots +B_{0}(u)F(u)=0,
\end{equation*}%
where $B_{0}(u)=\dfrac{1}{\alpha _{k}}A_{0}\left( z\left( u\right) \right) $
and 
\begin{equation*}
B_{j}(u)=\frac{\alpha _{j}}{\alpha _{k}}+\frac{\alpha _{j}}{\alpha _{k}}%
\sum_{n=j}^{k-1}A_{n}\left( z\left( u\right) \right) ,\text{ }j=1,2,\dots
,k-1.
\end{equation*}%
By the proof of Lemma 3.3, we can get that $\left[ 21,\text{ p. 63}\right] $%
\begin{equation*}
\alpha _{k}=V^{k}\left( u\right) =\left( \frac{1}{z^{\prime }(u)}\right) ^{k}
\end{equation*}%
\begin{equation*}
=\left( \frac{\omega }{e^{i\theta _{0}}}\left( \frac{1-u}{-\left( 1+u\right)
+\sqrt{2\left( 1+u^{2}\right) }}\right) ^{\frac{1}{\omega }-1}\frac{\left(
1-u\right) ^{2}\sqrt{1+u^{2}}}{\sqrt{2}\left( 1+u\right) -2\sqrt{1+u^{2}}}%
\right) ^{k},
\end{equation*}%
which is analytic in $\Delta ,$ where $\theta _{0}=\frac{\alpha +\beta }{2}$
and $\omega =\frac{\pi }{\beta -\alpha }.$ Since $\alpha _{k}=V^{k}\left(
u\right) \neq 0$ in $\Delta ,$ then $B_{0}(u)=\dfrac{1}{\alpha _{k}}%
A_{0}\left( z\left( u\right) \right) $ and 
\begin{equation*}
B_{j}(u)=\frac{\alpha _{j}}{\alpha _{k}}+\frac{\alpha _{j}}{\alpha _{k}}%
\sum_{n=j}^{k-1}A_{n}\left( z\left( u\right) \right) ,\text{ }j=1,2,\dots
,k-1
\end{equation*}%
are also analytic in $\Delta .$ Because 
\begin{equation*}
T\left( \varrho ,\alpha _{j}\right) =O\left( \log \frac{1}{1-\varrho }%
\right) ,\text{ \ }j=1,2,\dots ,k,
\end{equation*}%
it follows from this and the properties of Nevanlinna's characteristic
function that%
\begin{equation*}
T\left( \varrho ,B_{0}\right) \leq T\left( \varrho ,\frac{1}{\alpha _{k}}%
\right) +T\left( \varrho ,A_{0}\left( z\left( u\right) \right) \right)
\end{equation*}%
\begin{equation*}
=T\left( \varrho ,\alpha _{k}\right) +T\left( \varrho ,A_{0}\left( z\left(
u\right) \right) \right) +O\left( 1\right)
\end{equation*}%
\begin{equation*}
=T\left( \varrho ,A_{0}\left( z\left( u\right) \right) \right) +O\left( \log 
\frac{1}{1-\varrho }\right) ,
\end{equation*}%
and for $j=1,2,\dots ,k-1$%
\begin{equation*}
T\left( \varrho ,B_{j}\right) \leq T\left( \varrho ,\frac{\alpha _{j}}{%
\alpha _{k}}\right) +\sum_{n=j}^{k-1}T\left( \varrho ,A_{n}\left( z\left(
u\right) \right) \right) +O\left( 1\right)
\end{equation*}%
\begin{equation*}
\leq T\left( \varrho ,\alpha _{j}\right) +T\left( \varrho ,\frac{1}{\alpha
_{k}}\right) +\sum_{n=j}^{k-1}T\left( \varrho ,A_{n}\left( z\left( u\right)
\right) \right) +O\left( 1\right)
\end{equation*}%
\begin{equation*}
=T\left( \varrho ,\alpha _{j}\right) +T\left( \varrho ,\alpha _{k}\right)
+\sum_{n=j}^{k-1}T\left( \varrho ,A_{n}\left( z\left( u\right) \right)
\right) +O\left( 1\right)
\end{equation*}%
\begin{equation*}
=\sum_{n=j}^{k-1}T\left( \varrho ,A_{n}\left( z\left( u\right) \right)
\right) +O\left( \log \frac{1}{1-\varrho }\right) .
\end{equation*}%
\textbf{Lemma 3.5 }$\left[ 16\right] $\textbf{\ }\textit{Let }$p\geq q\geq 1$%
\textit{\ be integers. If }$B_{0}(u),B_{1}(u),...,B_{k-1}(u)$\textit{\ are
analytic functions of }[\textit{p, q}]\textit{-order in the unit disc }$%
\Delta $\textit{, then every solution }$F\not\equiv 0$\textit{\ of }$\left(
3.7\right) $\textit{\ satisfies}%
\begin{equation*}
\mu _{\lbrack p+1,q]}\left( F\right) =\mu _{M,[p+1,q]}\left( F\right) \leq
\max_{1\leq j\leq k-1}\left\{ \mu _{M,[p,q]}\left( B_{0}\right) ,\rho
_{M,[p,q]}\left( B_{j}\right) \right\} .
\end{equation*}%
\textbf{Lemma 3.6 }\textit{Let }$p\geq q\geq 1$\textit{\ be integers. If }$%
A_{0}(z),...,A_{k-1}(z)$\textit{\ are analytic functions of }[\textit{p, q}]%
\textit{-order in sector }$\Omega $\textit{\ satisfying }$\max\limits_{1\leq
j\leq k-1}\left\{ \rho _{\lbrack p,q],\Omega }\left( A_{j}\right) \right\}
<\mu _{\lbrack p,q],\Omega _{\varepsilon }}\left( A_{0}\right) ,$\textit{\
then for any given }$\varepsilon \in \left( 0,\frac{\beta -\alpha }{2}%
\right) ,$\textit{\ every solution }$f\not\equiv 0$\textit{\ of }$\left(
1.1\right) $\textit{\ satisfies}%
\begin{equation*}
\mu _{\lbrack p+1,q],\Omega _{\varepsilon }}\left( f\right) \leq \mu
_{\lbrack p,q],\Omega }\left( A_{0}\right) +1.
\end{equation*}%
\textit{Furthermore, if }$p>q$\textit{\ then}%
\begin{equation*}
\mu _{\lbrack p+1,q],\Omega _{\varepsilon }}\left( f\right) \leq \mu
_{\lbrack p,q],\Omega }\left( A_{0}\right) .
\end{equation*}%
\textit{Proof}. Let $f\not\equiv 0$ be a solution of equation $\left(
1.1\right) $. Then by Lemma 3.4, $\ F(u)=f\left( z\left( u\right) \right) $
is a solution of equation $\left( 3.7\right) $ and by using Remark 3.1,
Proposition 1.1, Proposition 1.2 and Lemma 3.5, we obtain%
\begin{equation*}
\mu _{\lbrack p+1,q],\Omega _{\varepsilon }}\left( f\right) \leq \mu
_{\lbrack p+1,q]}\left( F\right) =\mu _{M,[p+1,q]}\left( F\right)
\end{equation*}%
\begin{equation*}
\leq \max_{1\leq j\leq k-1}\left\{ \mu _{M,[p,q]}\left( B_{0}\right) ,\rho
_{M,[p,q]}\left( B_{j}\right) \right\}
\end{equation*}%
\begin{equation*}
\leq \max_{1\leq j\leq k-1}\left\{ \mu _{\lbrack p,q]}\left( B_{0}\right)
,\rho _{\lbrack p,q]}\left( B_{j}\right) \right\} +1
\end{equation*}%
\begin{equation*}
\leq \max_{1\leq j\leq k-1}\left\{ \mu _{\lbrack p,q],\Omega }\left(
A_{0}\right) ,\rho _{\lbrack p,q],\Omega }\left( A_{j}\right) \right\} +1
\end{equation*}%
\begin{equation*}
\leq \max_{1\leq j\leq k-1}\left\{ \mu _{\lbrack p,q],\Omega }\left(
A_{0}\right) ,\mu _{\lbrack p,q],\Omega _{\varepsilon }}\left( A_{0}\right)
\right\} +1=\mu _{\lbrack p,q],\Omega }\left( A_{0}\right) +1.
\end{equation*}%
If $\ p>q,$ we obtain%
\begin{equation*}
\mu _{\lbrack p+1,q],\Omega _{\varepsilon }}\left( f\right) \leq \mu
_{\lbrack p+1,q]}\left( F\right) =\mu _{M,[p+1,q]}\left( F\right)
\end{equation*}%
\begin{equation*}
\leq \max_{1\leq j\leq k-1}\left\{ \mu _{M,[p,q]}\left( B_{0}\right) ,\rho
_{M,[p,q]}\left( B_{j}\right) \right\}
\end{equation*}%
\begin{equation*}
=\max_{1\leq j\leq k-1}\left\{ \mu _{\lbrack p,q]}\left( B_{0}\right) ,\rho
_{\lbrack p,q]}\left( B_{j}\right) \right\}
\end{equation*}%
\begin{equation*}
\leq \max_{1\leq j\leq k-1}\left\{ \mu _{\lbrack p,q],\Omega }\left(
A_{0}\right) ,\rho _{\lbrack p,q],\Omega }\left( A_{j}\right) \right\}
\end{equation*}%
\begin{equation*}
\leq \max_{1\leq j\leq k-1}\left\{ \mu _{\lbrack p,q],\Omega }\left(
A_{0}\right) ,\mu _{\lbrack p,q],\Omega _{\varepsilon }}\left( A_{0}\right)
\right\} =\mu _{\lbrack p,q],\Omega }\left( A_{0}\right) .
\end{equation*}%
\textbf{Lemma 3.7 }$\left[ 7,15\right] $\textbf{\ }\textit{Let }$f$\textit{\
be a meromorphic function in the unit disc }$\Delta $\textit{\ and let }$%
k\in 
\mathbb{N}
$\textit{. Then }%
\begin{equation*}
m\left( r,\frac{f^{(k)}}{f}\right) =S(r,f),
\end{equation*}%
\textit{where }$S(r,f)=O\left( \log ^{+}T(r,f)+\log \left( \frac{1}{1-r}%
\right) \right) $\textit{, possibly outside a set }$F\subset \lbrack 0,1)$%
\textit{\ with }$\int_{F}\frac{dr}{1-r}<\infty $\textit{.}

\quad

\noindent \textbf{Lemma 3.8 }$\left[ 1,7\right] $\textbf{\ }\textit{Let }$%
g:\left( 0,1\right) \rightarrow \mathbf{%
\mathbb{R}
}$\textit{\ \ and }$h:\left( 0,1\right) \rightarrow \mathbf{%
\mathbb{R}
}$ \textit{be monotone increasing functions such that }$g\left( r\right)
\leq h\left( r\right) $\textit{\ holds outside of an exceptional set }$%
E\subset \lbrack 0,1)$ \textit{for which }$\int_{E}\frac{dr}{1-r}<\infty $%
\textit{.} \textit{Then} \textit{there exists} \textit{a constant} $d\in
\left( 0,1\right) $ \textit{such that if} $s\left( r\right) =1-d\left(
1-r\right) ,$ \textit{then} $g\left( r\right) \leq h\left( s\left( r\right)
\right) $\textit{\ for all }$r\in \lbrack 0,1).$

\quad

\noindent \textbf{Lemma 3.9 }$\left[ 25\right] $\textbf{\ }\textit{Let }$%
p\geq q\geq 1$\textit{\ be integers. If }$A_{0}(z),...,A_{k-1}(z)$\textit{\
are analytic functions of [p, q]-order in sector }$\Omega $\textit{\
satisfying }$\max\limits_{0\leq j\leq k-1}\left\{ \rho _{\lbrack p,q],\Omega
}\left( A_{j}\right) \right\} \leq \eta ,$\textit{\ then for any given }$%
\varepsilon \in \left( 0,\frac{\beta -\alpha }{2}\right) ,$\textit{\ every
solution }$f\not\equiv 0$\textit{\ of }$\left( 1.1\right) $\textit{\
satisfies}%
\begin{equation*}
\rho _{\lbrack p+1,q],\Omega _{\varepsilon }}\left( f\right) \leq \eta +1.
\end{equation*}%
\textit{Furthermore, if }$p>q$\textit{\ then}%
\begin{equation*}
\rho _{\lbrack p+1,q],\Omega _{\varepsilon }}\left( f\right) \leq \eta .
\end{equation*}

\section{\textbf{Proofs of the Theorems}}

\noindent \textbf{Proof of Theorem 2.1.} Suppose that $f\not\equiv 0$ is a
solution of $\left( 1.1\right) $ in the sector $\Omega .$ From Lemma 3.4,
the function $F\left( u\right) =f\left( z\left( (u\right) \right) $ is a
solution of $\left( 3.7\right) ,$ where $z\left( u\right) $ is defined by $%
\left( 3.2\right) .$ Then, by Lemma 3.2 and the properties of characteristic
function of Nevanlinna, we have%
\begin{equation*}
T\left( \varrho ,B_{0}(u)\right) =T\left( \varrho ,\frac{1}{\alpha _{k}}%
A_{0}\left( z\left( u\right) \right) \right) \geq T\left( \varrho
,A_{0}\left( z\left( u\right) \right) \right) -T\left( \varrho ,\alpha
_{k}\right) 
\end{equation*}%
\begin{equation*}
=T_{0}\left( \varrho ,\mathbb{C},A_{0}\left( z\left( u\right) \right)
\right) +O(1)-T\left( \varrho ,\alpha _{k}\right) 
\end{equation*}%
\begin{equation}
\geq \frac{b}{2}T_{0}\left( 1-\frac{1-\varrho }{b},\Omega _{\varepsilon
},A_{0}(z)\right) +O(1)-T\left( \varrho ,\alpha _{k}\right) .  \tag{4.1}
\end{equation}%
By $\left( 3.3\right) $, $\left( 3.11\right) $ and the formula $T\left(
r,f\right) =T_{0}\left( r,\mathbb{C},f\right) +O(1)$ $\left( 0<r<1\right) ,$
for $j=1,2,\dots ,k-1$ we have%
\begin{equation*}
T\left( \varrho ,B_{j}\left( u\right) \right) \leq \sum_{n=j}^{k-1}T\left(
\varrho ,A_{n}\left( z\left( u\right) \right) \right) +O\left( \log \frac{1}{%
1-\varrho }\right) 
\end{equation*}%
\begin{equation*}
=\sum_{n=j}^{k-1}T_{0}\left( \varrho ,\mathbb{C},A_{n}\left( z\left(
u\right) \right) \right) +O(1)+O\left( \log \frac{1}{1-\varrho }\right) 
\end{equation*}%
\begin{equation}
\leq \frac{16\pi }{\delta }\sum_{n=j}^{k-1}T_{0}\left( 1-\frac{\delta }{8\pi 
}\left( 1-\varrho \right) ,\Omega ,A_{n}\left( z\right) \right) +O\left(
\log \frac{1}{1-\varrho }\right) .  \tag{4.2}
\end{equation}%
Set%
\begin{equation*}
\eta =\max_{1\leq j\leq k-1}\left\{ \rho _{\lbrack p,q],\Omega }\left(
A_{j}\right) \right\} <\mu _{\lbrack p,q],\Omega _{\varepsilon }}\left(
A_{0}\right) =\mu .
\end{equation*}%
Then, for any given $\epsilon \left( 0<2\epsilon <\mu -\eta \right) $ and $%
r\rightarrow 1^{-},$ we have for $j=1,2,\dots ,k-1$ 
\begin{equation}
T_{0}\left( r,\Omega ,A_{j}(z)\right) \leq \exp _{p}\left\{ \left( \eta
+\epsilon \right) \log _{q}\frac{1}{1-r}\right\} .  \tag{4.3}
\end{equation}%
By the definition of lower $\left[ p,q\right] $ order 
\begin{equation}
T_{0}\left( r,\Omega _{\varepsilon },A_{0}(z)\right) \geq \exp _{p}\left\{
\left( \mu -\epsilon \right) \log _{q}\frac{1}{1-r}\right\} .  \tag{4.4}
\end{equation}%
Now, as $|u|=\varrho \rightarrow 1^{-},$ it follows from $\left( 4.1\right) ,
$ $\left( 4.2\right) ,$ $\left( 4.3\right) $ and $\left( 4.4\right) $ that%
\begin{equation*}
T\left( \varrho ,B_{0}\right) \geq \frac{b}{2}T_{0}\left( 1-\frac{1-\varrho 
}{b},\Omega _{\varepsilon },A_{0}(z)\right) +O(1)-T\left( \varrho ,\alpha
_{k}\right) 
\end{equation*}%
\begin{equation*}
\geq \frac{b}{2}\exp _{p}\left\{ \left( \mu -\epsilon \right) \log
_{q}\left( \frac{b}{1-\varrho }\right) \right\} +O(1)-T\left( \varrho
,\alpha _{k}\right) 
\end{equation*}%
\begin{equation}
=O\left( \exp _{p}\left\{ \left( \mu -\epsilon \right) \log _{q}\left( \frac{%
1}{1-\varrho }\right) \right\} \right) -T\left( \varrho ,\alpha _{k}\right) 
\tag{4.5}
\end{equation}%
and for $j=1,2,\dots ,k-1$%
\begin{equation*}
T\left( \varrho ,B_{j}\right) \leq \frac{16\pi }{\delta }\left( k-j\right)
\exp _{p}\left\{ \left( \eta +\epsilon \right) \log _{q}\left( \frac{8\pi }{%
\delta \left( 1-\varrho \right) }\right) \right\} +O\left( \log \frac{1}{%
1-\varrho }\right) 
\end{equation*}%
\begin{equation}
=O\left( \exp _{p}\left\{ \left( \eta +\epsilon \right) \log _{q}\left( 
\frac{1}{1-\varrho }\right) \right\} +\log \frac{1}{1-\varrho }\right) . 
\tag{4.6}
\end{equation}%
By $\left( 3.7\right) ,$ we can write%
\begin{equation*}
T\left( \varrho ,B_{0}\right) =m\left( \varrho ,B_{0}\right) \leq
\sum_{j=1}^{k-1}m\left( \varrho ,B_{j}\right) +\sum_{j=1}^{k}m\left( \varrho
,\frac{F^{(j)}}{F}\right) +O(1)
\end{equation*}%
\begin{equation}
=\sum_{j=1}^{k-1}T\left( \varrho ,B_{j}\right) +\sum_{j=1}^{k}m\left(
\varrho ,\frac{F^{(j)}}{F}\right) +O(1).  \tag{4.7}
\end{equation}%
It follows by $\left( 4.5\right) ,\left( 4.6\right) ,\left( 4.7\right) $ and
Lemma 3.7 that%
\begin{equation*}
O\left( \exp _{p}\left\{ \left( \mu -\epsilon \right) \log _{q}\left( \frac{1%
}{1-\varrho }\right) \right\} \right) \leq O\left( \exp _{p}\left\{ \left(
\eta +\epsilon \right) \log _{q}\left( \frac{1}{1-\varrho }\right) \right\}
\right) 
\end{equation*}%
\begin{equation}
+O\left( \log \frac{1}{1-\varrho }\right) +T\left( \varrho ,\alpha
_{k}\right) +O\left( \log ^{+}T\left( \varrho ,F\right) +\log \frac{1}{%
1-\varrho }\right)   \tag{4.8}
\end{equation}%
holds for all $u$ satisfying $|u|=\varrho \notin E$ as $\varrho \rightarrow
1^{-}$ and $E\subset (0,1)$ is a set with $\int_{E}\frac{d\varrho }{%
1-\varrho }<+\infty .$ By using Lemma 3.8 and $\left( 4.8\right) ,$ for all $%
u$ satisfying $|u|=\varrho $ as $\varrho \rightarrow 1^{-},$ we obtain%
\begin{equation*}
\exp _{p}\left\{ \left( \mu -\epsilon \right) \log _{q}\left( \frac{1}{%
1-\varrho }\right) \right\} \leq O\left( \exp _{p}\left\{ \left( \eta
+\epsilon \right) \log _{q}\left( \frac{1}{1-\varrho }\right) \right\}
\right) 
\end{equation*}%
\begin{equation}
+O\left( \log \frac{1}{d\left( 1-\varrho \right) }\right) +O\left( \log
^{+}T\left( 1-d\left( 1-\varrho \right) ,F\right) \right) .  \tag{4.9}
\end{equation}%
Thus, from $\left( 4.9\right) $ we get $\sigma _{\lbrack p,q]}\left(
F\right) =\mu _{\lbrack p+1,q]}\left( F\right) =+\infty $ and $\sigma
_{\lbrack p+1,q]}\left( F\right) \geq \mu _{\lbrack p+1,q]}\left( F\right)
\geq \mu .$ Then, by Remark 3.1, we get that 
\begin{equation*}
\rho _{\lbrack p,q],\Omega }\left( f\left( z\right) \right) =\mu _{\lbrack
p,q]}\left( f\left( z\right) \right) =+\infty \text{ and }\rho _{\lbrack
p+1,q],\Omega }\left( f\left( z\right) \right) \geq \mu _{\lbrack
p+1,q],\Omega }\left( f\left( z\right) \right) \geq \mu .
\end{equation*}%
On the other hand, by Lemma 3.6 we have $\mu _{\lbrack p+1,q],\Omega
_{\varepsilon }}\left( f\right) \leq \mu _{\lbrack p,q],\Omega }\left(
A_{0}\right) +1,$ and if $p>q,$ we have $\mu _{\lbrack p+1,q],\Omega
_{\varepsilon }}\left( f\right) \leq \mu _{\lbrack p,q],\Omega }\left(
A_{0}\right) .$

\quad

\noindent \textbf{Proof of Corollary 2.1. }By using Theorem 2.1 and Lemma
3.9, we easily obtain Corollary 2.1.

\quad

\noindent \textbf{Proof of Theorem 2.2.} Suppose that $f\not\equiv 0$ is a
solution of $\left( 1.1\right) $ in the sector $\Omega .$ From Lemma 3.4,
the function $F\left( u\right) =f\left( z\left( (u\right) \right) $ is a
solution of $\left( 3.7\right) ,$ where $z\left( u\right) $ is defined by $%
\left( 3.2\right) .$ If $\rho _{\lbrack p,q],\Omega }\left( A_{j}\right)
<\mu _{\lbrack p,q],\Omega _{\varepsilon }}\left( A_{0}\right) =\mu $ for
all $j=1,\cdots ,k-1$, then Theorem 2.2 reduces to Theorem 2.1. Thus, we
assume that at least one of $A_{j}$ $(j=1,\cdots ,k-1)$ satisfies $\rho
_{\lbrack p,q],\Omega }\left( A_{j}\right) =\mu _{\lbrack p,q],\Omega
_{\varepsilon }}\left( A_{0}\right) =\mu .$ So, there exists a set $%
I\subseteq \{1,\cdots ,k-1\}$ such that for $j\in I$ we have $\rho _{\lbrack
p,q],\Omega }\left( A_{j}\right) =\mu _{\lbrack p,q],\Omega _{\varepsilon
}}\left( A_{0}\right) =\mu $ and 
\begin{equation*}
\tau _{1}=\max_{j\in I}\{\tau _{\lbrack p,q],\Omega }\left( A_{j}\right)
:\rho _{\lbrack p,q],\Omega }\left( A_{j}\right) =\mu _{\lbrack p,q],\Omega
_{\varepsilon }}\left( A_{0}\right) \}<\underline{\tau }_{[p,q],\Omega
_{\varepsilon }}\left( A_{0}\right) =\tau <+\infty 
\end{equation*}%
and for $j\in \{1,\cdots ,k-1\}\backslash I,$ we have $b=\max_{j\in
\{1,\cdots ,k-1\}\backslash I}\{\rho _{\lbrack p,q],\Omega }\left(
A_{j}\right) \}<\mu _{\lbrack p,q],\Omega _{\varepsilon }}\left(
A_{0}\right) =\mu .$ Then that for any given $\epsilon $ $\left( 0<2\epsilon
<\min \left\{ \mu -b,\tau -\tau _{1}\right\} \right) $ and for $r\rightarrow
1^{-},$ we have for $j\in \{1,\cdots ,k-1\}\backslash I$ 
\begin{equation}
T_{0}\left( r,\Omega ,A_{j}(z)\right) \leq \exp _{p}\left\{ \left(
b+\epsilon \right) \log _{q}\frac{1}{1-r}\right\} \leq \exp _{p}\left\{
\left( \mu -\epsilon \right) \log _{q}\frac{1}{1-r}\right\}   \tag{4.10}
\end{equation}%
and for $j\in I,$ we get 
\begin{equation}
T_{0}\left( r,\Omega ,A_{j}(z)\right) \leq \exp _{p-1}\left\{ \left( \tau
_{1}+\epsilon \right) \left( \log _{q-1}\frac{1}{1-r}\right) ^{\mu }\right\}
.  \tag{4.11}
\end{equation}%
By the definition of lower $\left[ p,q\right] $ order, we have for $%
r\rightarrow 1^{-}$ 
\begin{equation}
T_{0}\left( r,\Omega _{\varepsilon },A_{0}(z)\right) \geq \exp _{p-1}\left\{
\left( \tau -\epsilon \right) \left( \log _{q-1}\frac{1}{1-r}\right) ^{\mu
}\right\} .  \tag{4.12}
\end{equation}%
Then, by $\left( 4.1\right) $ and $\left( 4.12\right) $ as $|u|=\varrho
\rightarrow 1^{-}$%
\begin{equation*}
T\left( \varrho ,B_{0}(u)\right) =T\left( \varrho ,\frac{1}{\alpha _{k}}%
A_{0}\left( z\left( u\right) \right) \right) 
\end{equation*}%
\begin{equation*}
\geq \frac{b}{2}T_{0}\left( 1-\frac{1-\varrho }{b},\Omega _{\varepsilon
},A_{0}(z)\right) +O(1)-T\left( \varrho ,\alpha _{k}\right) 
\end{equation*}%
\begin{equation*}
\geq \frac{b}{2}\exp _{p-1}\left\{ \left( \tau -\epsilon \right) \left( \log
_{q-1}\frac{b}{1-\varrho }\right) ^{\mu }\right\} +O(1)-T\left( \varrho
,\alpha _{k}\right) 
\end{equation*}%
\begin{equation}
=O\left( \exp _{p-1}\left\{ \left( \tau -\epsilon \right) \left( \log _{q-1}%
\frac{1}{1-\varrho }\right) ^{\mu }\right\} \right) -T\left( \varrho ,\alpha
_{k}\right) .  \tag{4.13}
\end{equation}%
Also, by $\left( 4.2\right) ,$ $\left( 4.10\right) $ and $\left( 4.11\right) 
$ for $j=1,2,\dots ,k-1$%
\begin{equation*}
T\left( \varrho ,B_{j}\right) \leq \frac{16\pi }{\delta }%
\sum_{n=j}^{k-1}T_{0}\left( 1-\frac{\delta }{8\pi }\left( 1-\varrho \right)
,\Omega ,A_{n}\left( z\right) \right) +O\left( \log \frac{1}{1-\varrho }%
\right) 
\end{equation*}%
\begin{equation*}
\leq O\left( \exp _{p}\left\{ \left( \mu -\epsilon \right) \log _{q}\frac{%
8\pi }{\delta \left( 1-\varrho \right) }\right\} \right) 
\end{equation*}%
\begin{equation*}
+O\left( \exp _{p-1}\left\{ \left( \tau _{1}+\epsilon \right) \left( \log
_{q-1}\frac{8\pi }{\delta \left( 1-\varrho \right) }\right) ^{\mu }\right\}
\right) +O\left( \log \frac{1}{1-\varrho }\right) 
\end{equation*}%
\begin{equation}
=O\left( \exp _{p-1}\left\{ \left( \tau _{1}+\epsilon \right) \left( \log
_{q-1}\frac{1}{1-\varrho }\right) ^{\mu }\right\} +\log \frac{1}{1-\varrho }%
\right) .  \tag{4.14}
\end{equation}%
It follows by $\left( 4.7\right) ,\left( 4.13\right) ,$ $\left( 4.14\right) $
and Lemma 3.7 that%
\begin{equation*}
O\left( \exp _{p-1}\left\{ \left( \tau -\epsilon \right) \left( \log _{q-1}%
\frac{1}{1-\varrho }\right) ^{\mu }\right\} \right) \leq O\left( \exp
_{p-1}\left\{ \left( \tau _{1}+\epsilon \right) \left( \log _{q-1}\frac{1}{%
1-\varrho }\right) ^{\mu }\right\} \right) 
\end{equation*}%
\begin{equation}
+O\left( \log \frac{1}{1-\varrho }\right) +T\left( \varrho ,\alpha
_{k}\right) +O\left( \log ^{+}T\left( \varrho ,F\right) +\log \frac{1}{%
1-\varrho }\right)   \tag{4.15}
\end{equation}%
holds for all $u$ satisfying $|u|=\varrho \notin E$ as $\varrho \rightarrow
1^{-},$ where $E\subset (0,1)$ is a set with $\int_{E}\frac{d\varrho }{%
1-\varrho }<+\infty .$ By using Lemma 3.8 and $\left( 4.15\right) ,$ for all 
$u$ satisfying $|u|=\varrho \rightarrow 1^{-},$ we obtain%
\begin{equation*}
\exp _{p-1}\left\{ \left( \tau -\epsilon \right) \left( \log _{q-1}\frac{1}{%
1-\varrho }\right) ^{\mu }\right\} \leq O\left( \exp _{p-1}\left\{ \left(
\tau _{1}+\epsilon \right) \left( \log _{q-1}\frac{1}{d\left( 1-\varrho
\right) }\right) ^{\mu }\right\} \right) 
\end{equation*}%
\begin{equation}
+O\left( \log \frac{1}{d\left( 1-\varrho \right) }\right) +O\left( \log
^{+}T\left( 1-d\left( 1-\varrho \right) ,F\right) \right) .  \tag{4.16}
\end{equation}%
Thus, from $\left( 4.16\right) $ we get $\rho _{\lbrack p,q]}\left( F\right)
=\mu _{\lbrack p,q]}\left( F\right) =+\infty $ and $\rho _{\lbrack
p+1,q]}\left( F\right) \geq \mu _{\lbrack p+1,q]}\left( F\right) \geq \mu .$
Then, by Remark 3.1, we get that 
\begin{equation*}
\rho _{\lbrack p,q],\Omega }\left( f\left( z\right) \right) =\mu _{\lbrack
p,q]}\left( f\left( z\right) \right) =+\infty \text{ and }\rho _{\lbrack
p+1,q],\Omega }\left( f\left( z\right) \right) \geq \mu _{\lbrack
p+1,q],\Omega }\left( f\left( z\right) \right) \geq \mu .
\end{equation*}%
On the other hand, by Lemma 3.6 we have $\mu _{\lbrack p+1,q],\Omega
_{\varepsilon }}\left( f\right) \leq \mu _{\lbrack p,q],\Omega }\left(
A_{0}\right) +1,$ and if $p>q,$ we have $\mu _{\lbrack p+1,q],\Omega
_{\varepsilon }}\left( f\right) \leq \mu _{\lbrack p,q],\Omega }\left(
A_{0}\right) .$

\quad

\noindent \textbf{Acknowledgements.} The author would like to thank the
anonymous referee and editor for their helpful remarks and suggestions to
improve this article. This paper is supported by University of Mostaganem
(UMAB) (PRFU Project Code C00L03UN270120180005).

\begin{center}
{\large References}
\end{center}

\noindent $\left[ 1\right] $ S. Bank, \textit{General theorem concerning the
growth of solutions of first-order algebraic differential equations},
Compositio Math. 25 (1972), 61--70.

\noindent $\left[ 2\right] $ B. Bela\"{\i}di, \textit{Growth of solutions to
linear equations with analytic coefficients of }[\textit{p, q}]\textit{%
-order in the unit disc}, Electron. J. Differential Equations 2011, No. 156,
1-11.

\noindent $\left[ 3\right] $ B. Bela\"{\i}di, \textit{On the }[\textit{p,q}]%
\textit{-order of analytic solutions of linear differential equations in the
unit disc}, Novi Sad J. Math. 42 (2012), no. 1, 117--129.

\noindent $\left[ 4\right] $ I. Chyzhykov, G. Gundersen, J. Heittokangas, 
\textit{\ Linear differential equations and logarithmic derivative estimates}%
, Proc. London Math. Soc. (3) 86 (2003), no. 3, 735--754.

\noindent $\left[ 5\right] $ A.A. Goldberg and I.V. Ostrovskii, \textit{%
Value distribution of meromorphic functions}, Translations of Mathematical
Monographs, 236. American Mathematical Society, Providence, RI, 2008.

\noindent $\left[ 6\right] $ W. K. Hayman, \textit{Meromorphic functions},
Oxford Mathematical Monographs Clarendon Press, Oxford, 1964.

\noindent $\left[ 7\right] $ J. Heittokangas, \textit{On complex
differential equations in the unit disc}, Ann. Acad. Sci. Fenn. Math. Diss.
No. 122 (2000), 54 pp.

\noindent $\left[ 8\right] $ H. Hu and X. M. Zheng, \textit{Growth of
solutions of linear differential equations with analytic coefficients of }[%
\textit{p,q}]\textit{-order in the unit disc}, Electron. J. Differential
Equations 2014, No. 204, 12 pp.

\noindent $\left[ 9\right] $ I. Laine, \textit{Complex differential equations%
}, Handbook of Differential Equations : Ordinary Differential Equations,
Vol. IV, 269--363, Handb. Differ. Equ., Elsevier/North-Holland, Amsterdam,
2008.

\noindent $\left[ 10\right] $ Z. Latreuch and B. Bela\"{\i}di, \textit{%
Linear differential equations with analytic coefficients of }[\textit{p,q}]%
\textit{-order in the unit disc}, Sarajevo J. Math. 9(21) (2013), no. 1,
71--84.

\noindent $\left[ 11\right] $ J. Long, \textit{Growth of solutions of higher
order complex linear differential equations in an angular domain of unit disc%
}, J. Math. Study 48 (2015), no. 3, 306--314.

\noindent $\left[ 12\right] $ J. Long, \textit{On} [\textit{p, q}]-\textit{%
order of solutions of higher-order complex linear differential equations in
an angular domain of unit disc}, J. Math. Study 50 (2017), no. 1, 91--100.

\noindent $\left[ 13\right] $ Ch. Pommerenke, \emph{\ }\textit{On the mean
growth of the solutions of complex linear differential equations in the disk}%
, Complex Variables Theory Appl. 1 (1982/83), no. 1, 23--38.

\noindent $\left[ 14\right] $ D. C. Sun and J. R. Yu, \textit{On the
distribution of random Dirichlet series}\emph{\ }$($II$)$, Chin. Ann. Math.
11(B) (1990), 33--44.

\noindent $\left[ 15\right] $ M. Tsuji, \textit{Potential Theory in Modern
Function Theory}, Chelsea, New York, (1975), reprint of the 1959 edition.

\noindent $\left[ 16\right] $ J. Tu and H. X. Huang, \textit{Complex
oscillation of linear differential equations with analytic coefficients of }[%
\textit{p,q}]\textit{-order in the unit disc}, Comput. Methods Funct. Theory
15 (2015), no. 2, 225--246.

\noindent $\left[ 17\right] $ S. J. Wu, \textit{Estimates for the
logarithmic derivative of a meromorphic function in an angle, and their
application}, Proceedings of the Conference on Complex Analysis (Tianjin,
1992), 235--240, Conf. Proc. Lecture Notes Anal., I, Int. Press, Cambridge,
MA, 1994.

\noindent $\left[ 18\right] $ S. J. Wu, \textit{On the growth of solutions
of second order linear differential equations in an angle}, Complex
Variables Theory Appl. 24 (1994), no. 3-4, 241--248.

\noindent $\left[ 19\right] $ N. Wu, \textit{Growth of solutions to linear
complex differential equations in an angular regio}n, Electron. J.
Differential Equations 2013, No. 183, 8 pp.

\noindent $\left[ 20\right] $ N. Wu and Y. Z. Li, \textit{On the growth of
solutions of higher order linear differential equations}, New Zealand J.
Math. 42 (2012), 27--35.

\noindent $\left[ 21\right] $ N. Wu, \textit{On the growth order of
solutions of linear differential equations in a sector of the unit disk},
Results Math. 65 (2014), no. 1-2, 57--66.

\noindent $\left[ 22\right] $ J. F. Xu and H. X. Yi, \textit{Solutions of
higher order linear differential equations in an angle}, Appl. Math. Lett.
22 (2009), no. 4, 484--489.

\noindent $\left[ 23\right] $ C. C. Yang and H. X. Yi, \textit{Uniqueness
theory of meromorphic functions}, Mathematics and its Applications, 557.
Kluwer Academic Publishers Group, Dordrecht, 2003.

\noindent $\left[ 24\right] $ G. Zhang, \textit{Value distributions of
solutions to complex linear differential equations in angular domains}, Open
Math. 15 (2017), no. 1, 884--894.

\noindent $\left[ 25\right] $ M. A. Zemirni and B. Bela\"{\i}di, [\textit{p,q%
}]\textit{-order of solutions of complex differential equations in a sector
of the unit disc}, An. Univ. Craiova Ser. Mat. Inform. 45 (2018), no. 1,
37--49.

\end{document}